\documentclass[twoside,12pt]{article}

\usepackage{latexsym,times,amssymb,mathptm}
\usepackage[all]{xy}

\def\demo{\noindent{\bf Proof. }}
\def\sqr#1#2{{\vcenter{\hrule height.#2pt
    \hbox{\vrule width.#2pt height#1pt \kern#1pt
        \vrule width.#2pt}
    \hrule height.#2pt}}}
\def\square{\mathchoice\sqr64\sqr64\sqr{4}3\sqr{3}3}
\def\QED{\hfill$\square$}

\def\tratto{\mbox{\rule{2mm}{.2mm}$\;\!$}}

\def\m{{\mathfrak m}}

\def\frak{\mathfrak}

\newfont{\frank}{eufm10 scaled\magstep1}

\newfont{\frankk}{eufm10 scaled\magstephalf}

\newfont{\bella}{eusm10 scaled\magstep1}

\def\II{\mbox{\bella \symbol{'111}}}
\def\JJ{\mbox{\bella \symbol{'112}}}

\def\MM{\mbox{\bella \symbol{'115}}}

\newtheorem{Theorem}{Theorem}[section]
\newtheorem{Lemma}[Theorem]{Lemma}
\newtheorem{Corollary}[Theorem]{Corollary}
\newtheorem{Proposition}[Theorem]{Proposition}

\newtheorem{Remark}[Theorem]{Remark}
\newtheorem{Example}[Theorem]{Example}

\newtheorem{Question}[Theorem]{Question}

\begin{document}

\title{
{\bf Sally Modules and Associated Graded Rings}}

\author{
\begin{tabular}{ccc}
{\bf Alberto Corso}\footnote{The first two authors were partially
supported by the NATO/CNR Advanced Fellowships Programme.}
& \ &
{\bf Claudia Polini}${}^{*}$ \\
\mbox{\normalsize Department of Mathematics}
& \ &
\mbox{\normalsize Department of Mathematics} \\
\mbox{\normalsize Purdue University}
& \ &
\mbox{\normalsize Michigan State University} \\
\mbox{\normalsize West Lafayette, IN 47907}
& \ &
\mbox{\normalsize East Lansing, MI 48824}  \\
\mbox{\normalsize E-mail: corso@math.purdue.edu}
& \ &
\mbox{\normalsize E-mail: polini@math.msu.edu}
\end{tabular}
\\ \ \\
\begin{tabular}{c}
{\bf Maria Vaz Pinto} \\
\mbox{\normalsize Centro de Algebra da Universidade de Lisboa} \\
\mbox{\normalsize Avenida Professor Gama Pinto, 2} \\
\mbox{\normalsize 1600 Lisboa, Portugal} \\
\mbox{\normalsize E-mail: vazpinto@hermite.cii.fc.ul.pt}
\end{tabular}}

\date{\ }

\maketitle

\vspace{-.5in}

\pagestyle{myheadings}

\markright{Sally modules and associated graded rings}
\markboth{A. Corso, C. Polini, and M. Vaz Pinto}{Sally modules and
associated graded rings}

\baselineskip 15 pt

\section{Introduction}

To frame and motivate the goals pursued in the present article we
recall that, loosely speaking, the most common among the blowup algebras
are the {\it Rees algebra} $R[It]=\bigoplus_{n=0}^{\infty} I^nt^n$ and the
{\it associated graded ring} $\mbox{\rm gr}_I(R)=\bigoplus_{n=0}^{\infty}
I^n/I^{n+1}$ of an ideal $I$ in a commutative Noetherian local ring
$(R, \m)$. The three main clusters around which most of the current
research on blowup algebras has been developed are: $(${\it a}$)$
the study of the depth properties of $R[It]$, or of an appropriate
object related to it such as its Proj;
$(${\it b}$)$ the comparison between the arithmetical properties of
$R[It]$ and $\mbox{\rm gr}_I(R)$;
$(${\it c}$)$ the correspondence between the Hilbert$/$Hilbert--Samuel
functions and the properties of $\mbox{\rm gr}_I(R)$ for an $\m$-primary
ideal $I$.

In this paper we address the relation mentioned in $(${\it c}$)$.
To make the terminology more precise, the
Hilbert--Samuel function is the numerical function that measures the growth
of the length of $R/I^n$, $\lambda(R/I^n)$, for all $n \geq 1$. It is
well known that for $n \gg 0$ this function is a polynomial
in $n$ of degree $d$, namely
\[
e_0 {n+d-1 \choose d} - e_1 {n+d-2 \choose d-1} + \cdots + (-1)^d e_d,
\]
where $d$ denotes the dimension of the ring $R$ and $e_0, e_1, \ldots, e_d$
are the normalized coefficients of the {\it Hilbert--Samuel
polynomial} of $I$.

Pioneering work on the interplay described in $(${\it c}$)$ was made by
Judith Sally in a sequence of papers \cite{S, S1, S2, S3, S3.5, S3.75,
S4}. A major
recognition of her important contribution came with the introduction of
the Sally module $($see \cite{Vas, Vas2}$)$.
Its definition requires a notion that has proved to be
quite useful in the theory of Rees algebras:
We say that the ideal $J\subset I$ is a {\it reduction} of $I$ if
$I^{n+1}=JI^n$ for some $n \geq 0$. The least such $n$ is called the
reduction number of $I$ with respect to $J$ and denoted by $r_J(I)$.
A {\it minimal reduction} is a reduction which is minimal with respect
to inclusion. Minimal reductions always exist and if the residue field
of the ring $R$ is infinite the number of generators of any minimal
reduction of $I$ equals the analytic spread of the ideal $I$, namely
$\mbox{\rm dim}(R[It] \otimes R/\m)$. If $I$ is an $\m$-primary ideal with a
minimal reduction $J$, the {\it Sally module of $I$ with respect to
$J$}, $S_J(I)$, is the graded $R[Jt]$-module, of dimension $d$ whenever
$S_J(I)\not= 0$, defined by the short exact sequence
\begin{equation}\label{sally}
0 \rightarrow I R[Jt] \longrightarrow I R[It] \longrightarrow
S_J(I) = \bigoplus_{n=2}^{\infty} I^n/J^{n-1}I \rightarrow 0.
\end{equation}
This new object is the outgrowth of a successful attempt made by W.V.
Vasconcelos to give a unified and at the same time simplified
treatment of several results by J. Sally and others. This comes about
as follows.
For $n \gg 0$, the growth of the length of the graded pieces of the Sally
module $S_J(I)$ is also measured by a polynomial in $n$ of degree
$d-1$
\[
s_0 {n+d-2 \choose d-1} - s_1 {n+d-3 \choose d-2} + \cdots +
(-1)^{d-1} s_{d-1}.
\]
The $e_i$'s relate to the $s_i$'s in the following manner $($see \cite{Vas,
Vas2}$)$
\[
e_0 = \lambda(R/J) \qquad e_1 = \lambda(I/J) + s_0 \qquad
e_i = s_{i-1} \quad {\rm for \ } i=2, \ldots, d
\]
so that: $(${\it a}$)$ $e_0 - e_1
\leq \lambda(R/I)$ $($due to Northcott, see \cite{No}$)$; $(${\it b}$)$
$e_0 - e_1 = \lambda(R/I)$ if and only if $I^2=JI$ $($due to
Huneke, see \cite{H}, and Ooishi, see \cite{O}$)$; $(${\it c}$)$ $e_1
\geq 0$. In this spirit, one of the results we give in Section 2 is a
simple proof of the positivity of $e_2$ $($due to Narita, see
\cite{Nar}$)$ and another lower bound for $e_1$ $($see
Proposition~\ref{Narita}$)$. We also show the independence from
the minimal reduction of the length of the graded components of the
Sally module $($see Proposition~\ref{independence}$)$.

A recurring theme in the work of J. Sally is the discovery of conditions
on the multiplicity $e$ of the local ring $(R, {\frak m})$ that assure
that $\mbox{\rm gr}_{\frak m}(R)$ is Cohen--Macaulay. By \cite{Aby}, the
multiplicity $e$ of $R$ satisfies the
inequality $e \geq \mu(\m)-d+1$, where $\mu(\m)$ denotes the minimal
number of generators of $\m$. More precisely, the closed formula is:
$e=\mu(\m)-d+1+\lambda(\m^2/J\m)$, where $J$ is a minimal reduction of $\m$.
If $\lambda(\m^2/J\m)=0$, i.e., $R$ has minimal multiplicity, J. Sally
proved in \cite{S} that $\mbox{\rm gr}_{\frak m}(R)$ is always Cohen--Macaulay.
After this case was settled it was natural to investigate the case in
which $\lambda({\frak m}^2/J{\frak m})=1$, i.e., $e = \mu({\frak m})-d+2$.
In \cite{S3.5} she proved that if in addition $R$ is Gorenstein then
$\mbox{\rm gr}_{\frak m}(R)$ is also Gorenstein.  Later, in
\cite{S4} she established the Cohen--Macaulayness of $\mbox{\rm gr}_{\frak
m}(R)$ for an arbitrary Cohen--Macaulay ring $R$ having type $s$
different from $\mu({\frak m})-d$ $($we recall that the type of a
Cohen--Macaulay local ring $(R, {\frak m})$ of dimension $d$ is given
by $\mbox{\rm dim}_{R/{\frak m}}(\mbox{\rm Ext}_R^d(R/{\frak m}, R)))$.
Still in \cite{S4} she exhibited examples of rings having type
$\mu({\frak m})-d$ and with $\mbox{\rm gr}_{\frak m}(R)$ not
Cohen--Macaulay; however, in all the given examples $\mbox{\rm
depth}(\mbox{\rm gr}_{\frak m}(R))$ always turned out to be $d-1$.
Therefore the conjecture that arose from this kind of scenario was
whether, in the critical case, the depth of $\mbox{\rm gr}_{\frak m}(R)$
is always at least $d-1$.
A simpler proof of Sally's results was given in the Ph.D.
theses of two of the authors $($see \cite{Vaz, P}$)$. There, they also
verified the conjecture with the additional assumption that the reduction
number of ${\frak m}$ with respect to $J$ is at most $4$.
Finally, in 1996, M.E. Rossi and G. Valla $($see \cite{RV}$)$ and  H. Wang
$($see \cite{W2}$)$ positively solved, at the same time, Sally's conjecture
using two different methods. Based on the proof of Rossi--Valla, later
S. Huckaba proved that if $\lambda(\m^3/J\m^2) \leq 1$ then $\mbox{\rm
gr}_{\frak m}(R)$ has depth at least $d-1$. In fact, he showed that
the same conclusion holds for any $\m$-primary ideal $I$ such that $J
\cap I^2=JI$ and $\lambda(I^3/JI^2) \leq 1$.

The original trust of our work was to see to which extent one could
generalize the above results. The main theorems of the paper
appear in the third section and, roughly speaking, deal with the class
of ${\mathfrak m}$-primary ideals $I$ in a Cohen--Macaulay $($sometimes even
Gorenstein$)$ local ring $(R, {\mathfrak m})$ such that: $(${\it a}$)$
$J \cap I^k=JI^{k-1}$ for $k=1, \ldots, n$; and $(${\it b}$)$
$\lambda(I^{n+1}/JI^n) = 1$. To be more precise, we show:

\medskip

\noindent{\bf Theorem~\ref{usingtheform}}
{\it Let $(R, \m)$ be a Cohen--Macaulay local ring of dimension $d>0$ and
infinite residue field. Let $I$ be an $\m$-primary ideal of $R$
and let $J$ be a minimal reduction of $I$ with
$\lambda(I^{n+1}/JI^n)=1$ for some $n\geq 1$. If the following hold
\begin{itemize}
\item[$(${\it a}$)$]
$J \cap I^k = JI^{k-1}$ for all $k=1, \ldots, n$;

\item[$(${\it b}$)$]
the vector space dimension of $V=I+J\colon I^n/J\colon I^n$ is at least $2$;
\end{itemize}
then the associated graded ring of $I$ is Cohen--Macaulay. }

\medskip

\noindent{\bf Theorem~\ref{secondmain}}
{\it Let $(R, {\mathfrak m})$ be a Gorenstein local ring of dimension $d>0$
and infinite residue field. Let $J$ be a minimal reduction of
${\mathfrak m}$ with $\lambda({\mathfrak m}^3/J{\mathfrak m}^2)=1$.
Then the associated graded ring of ${\mathfrak m}$ is
Cohen--Macaulay. }

\medskip

\noindent{\bf Theorem~\ref{main}}
{\it Let $(R, {\frak m})$ be a Cohen--Macaulay local ring of dimension
$d > 0$ and infinite residue field. Let $I$ be an ${\frak m}$-primary
ideal of $R$ with a minimal reduction $J$ satisfying the following conditions
\begin{itemize}
\item[$(${\it a}$)$]
$J \cap I^k = JI^{k-1}$ for all $k=1, \ldots, n$;

\item[$(${\it b}$)$]
$\lambda(I^{n+1}/JI^n) \leq 1$.
\end{itemize}
Then $\mbox{\rm depth}(\mbox{\rm gr}_I(R)) \geq d-1$. }

\medskip

The results stated in the previous theorems require conditions on the
length of $I^{n+1}/JI^n$, where $I$ is an ${\mathfrak m}$-primary
ideal with minimal reduction $J$. It is therefore natural to investigate the
independence, from the minimal reduction, of such lengths.
We show in Proposition~\ref{higher-lengths} that the depth of
$\mbox{\rm gr}_I(R)$ being at least $d-1$ is a sufficient
condition. The independence of these lengths was first observed by
T. Marley in \cite[Corollary 2.9]{M0} and then later recovered by
S. Huckaba in \cite[Corollary 2.6]{Huc1}. In each case, the means of
the proof are different: Our proof is a consequence of the
existence and the properties of a natural filtration of the Sally
module $S_J(I)$ introduced in \cite{Vaz, Vaz2}. The result about the
independence of $\lambda(R/J\cap I^n)$ is instead still open.

Section 4 ends the paper by describing various classes of ideals
where the hypotheses required in Theorem~\ref{main} are satisfied.

\medskip

Throughout the paper, the notation and terminology are the ones of
\cite{BH} and \cite{Mat}.

\subsection*{Acknowledgements}

The authors sincerely thank Luisa Doering, Bernd Ulrich, and Wolmer
V. Vasconcelos for useful discussions they had concerning the material
in this paper.

\section{On the Hilbert--Samuel coefficients, Sally modules, and
independence of lengths}

Given any minimal reduction $J$ of an $\m$-primary ideal $I$ one can
define, using $($\ref{sally}$)$, the Sally module of $I$ with respect
of $J$, $S_J(I)$. The following proposition shows, in particular, that
the length of each graded component is in fact an invariant of $I$.

\begin{Proposition}\label{independence}
Let $(R, \m)$ be a Cohen--Macaulay local ring of dimension $d>0$ and
infinite residue field $R/\m$. Let $I$ be an $\m$-primary ideal of $R$
and let $J$ be a minimal reduction of $I$. Then the following lengths
are independent of $J$
\begin{itemize}
\item[$(${\it a}$)$]
$\lambda(I/J)$ and $\lambda(I^n/J^{n-1}I)$ for $n\geq 2$;

\item[$(${\it b}$)$]
$\lambda(R/J\colon I)$;

\item[$(${\it c}$)$]
$\lambda(S_2(I/J))$, where $S_2(I/J)$ denotes the second component of
the symmetric algebra of $I/J$.
\end{itemize}
\end{Proposition}
\demo
$(${\it a}$)$
We can write $\lambda(I/J)=\lambda(R/J)-\lambda(R/I)=e_0-\lambda(R/I)$,
where $e_0$ denotes the multiplicity of $I$.
For $n \geq 2$ one has that $I^n/J^{n-1}I$ is the component of degree $n-1$
of the Sally module of $I$ with respect to $J$. From \cite{Vas,
Vas2, Vaz} it follows that
\[
\lambda(I^n/J^{n-1}I) = e_0 {n+d-2 \choose d} + \lambda(R/I)
{n+d-2 \choose d} - \lambda(R/I^n).
\]

$(${\it b}$)$
The claim follows from the fact that $\lambda(R/J\colon I) =
\lambda(R/J) - \lambda(J\colon I/J) = e_0 - \lambda(H_{m-d}(I))$,
where $m$ denotes the minimal number of generators of $I$ and
$H_{m-d}(I)$ is the last nonvanishing Koszul homology of $I$.

$(${\it c}$)$
One has the following short exact sequence
\[
0 \rightarrow \delta(I) \longrightarrow S_2(I/J) \longrightarrow
I^2/JI \rightarrow 0,
\]
where $\delta(I)$ denotes the kernel of the natural surjection from
$S_2(I)$ to $I^2$ $($see \cite[Proof of Remark 2.7]{CP}$)$.
Now again the assertion follows from $(${\it a}$)$.
\QED

\medskip

\begin{Remark}
{\rm Proposition~\ref{independence} shows that $\lambda(I/J)$ and
$\lambda(I^2/JI)$ never depends on the minimal reduction $J$ of $I$.
It is natural then to address the issue of the independence from $J$
of $\lambda(I^n/JI^{n-1})$ for
any $n$. In general, though, such an independence fails for $n \geq 3$
as the reduction number $r$ of $I$ may depend on the chosen minimal
reduction $J$, unless the depth of the associated graded ring of
$I$ is at least $d-1$ $($see \cite{Huc, M, T}$)$.
On the other hand, the sum of the previous lengths gives an upper
bound for the coefficient $e_1$ of the Hilbert--Samuel polynomial
of $I$. To be more precise one has the inequalities
\begin{equation}\label{inequalities}
\sum_{n=1}^{r} \lambda((I^n,J)/J) \leq e_1 \leq \sum_{n=1}^r
\lambda(I^n/JI^{n-1})
\end{equation}
$($see \cite{HM, Huc1, Vaz}$)$ and \cite{Huc1} shows that
$e_1=\sum_{n=1}^r \lambda(I^n/JI^{n-1})$ if and only if the depth of the
associated graded ring of $I$ is at least $d-1$.
In Proposition~\ref{higher-lengths} we show that this lower
bound on the depth of the associated graded ring of $I$ is also
sufficient for the independence of each single length and not just of the
entire sum $($see also  \cite[Corollary 2.9]{M0} and \cite[Corollary
2.6]{Huc1}$)$. }
\end{Remark}

\begin{Proposition}\label{higher-lengths}
Let $(R, \m)$ be a Cohen--Macaulay local ring of dimension $d>0$ and
infinite residue field $R/\m$. Let $I$ be an $\m$-primary ideal of $R$
and let $J$ be a minimal reduction of $I$. If
$\mbox{\rm depth}(\mbox{\rm gr}_I(R)) \geq d-1$ then
$\lambda(I^n/JI^{n-1})$ does not depend on $J$, for all $n \geq 1$.
\end{Proposition}
\demo
By \cite{Vaz2}, $\mbox{depth}(\mbox{gr}_I(R)) \geq d-1$ is
equivalent to $s_0 = \sum_{n=2}^r
\lambda(I^n/JI^{n-1})$, where $s_0$ is the multiplicity of $S_J(I)$.
Hence the Hilbert--Poincar\'e series of $\mbox{gr}_I(R)$ has the form
\begin{equation}\label{hilb-poin}
HP(\mbox{gr}_I(R), t) =
\frac{\lambda(R/I) + \sum\limits_{n=1}^r \left[\lambda(I^n/JI^{n-1}) -
\lambda(I^{n+1}/JI^n) \right] t^n}{(1-t)^d}
\end{equation}
$($see \cite[Corollary 1.2]{Vaz2}$)$. In particular, the numerator of
$($\ref{hilb-poin}$)$ is a polynomial, say $p(t)$, with coefficients
independent from $J$.

We now proceed by induction on $n$. The cases $n=0, 1$ are
taken care of by Proposition~\ref{independence}. For $n \geq 2$
consider the identity
\[
\lambda(I^{n+1}/JI^n) = \lambda(I^n/JI^{n-1}) - \left[
\lambda(I^n/JI^{n-1}) - \lambda(I^{n+1}/JI^n) \right];
\]
the assertion now follows by induction. \QED

\medskip

\begin{Remark}
{\rm The first inequality in $($\ref{inequalities}$)$ also raises the
issue of whether or not the condition $\mbox{depth}(\mbox{gr}_I(R)) \geq d-1$
guarantees that $\lambda((I^n,J)/J) = \lambda(I^n/J\cap I^n)$ does not
depend on $J$ for $n \geq 0$. By Proposition~\ref{higher-lengths},
the question has a positive answer if the associated graded ring
of $I$ is Cohen--Macaulay $($see \cite[Corollary 2.7]{VV}$)$. The case
in which $\mbox{depth}(\mbox{gr}_I(R))=d-1$ is still open. }
\end{Remark}

\smallskip

We end this section by proving two results on the normalized
Hilbert--Samuel coefficients of $I$ by means of the Sally module. The
first result gives a lower bound for $e_1$; it is in general less
sharp than the one in $($\ref{inequalities}$)$ but it has the
advantage of being independent of the minimal reduction of $I$. The
second one is a simpler proof of the positivity of $e_2$.

\begin{Proposition}\label{Narita}
Let $(R, \m)$ be a Cohen--Macaulay local ring of dimension $d>0$. Let
$I$ be an $\m$-primary ideal of $R$. Then the following hold
\begin{itemize}
\item[$(${\it a}$)$]
$e_1 \geq 2 e_0 - \lambda(R/I^2)$;

\item[$(${\it b}$)$]
$(${\rm Narita, see \cite{Nar}}$)$ if $d\geq 2$ then $e_2 \geq 0$. Moreover, if
$d=2$ then $e_2=0$ if and only if $I^q$ has reduction number one
for $q\gg 0$.
\end{itemize}
\end{Proposition}
\demo
$(${\it a}$)$
By passing to the faithfully flat extension $R \longrightarrow
R[X]_{{\frak m}[X]}$ we may assume that the residue field of $R$ is
infinite and that $J$ is a minimal reduction of $I$.
After tensoring with $R/I$ the defining sequence $($\ref{sally}$)$
of $S_J(I)$ one has
\[
I/I^2 R[Jt] \longrightarrow \mbox{gr}_I(R)_{+} = \bigoplus_{i =
1}^{\infty} I^i/I^{i+1} \longrightarrow S_J(I) \otimes R/I \rightarrow 0.
\]
Notice that the dimension of $S_J(I) \otimes R/I$ is either $d$ or
$0$, since the set of associated primes of $S_J(I)$, as
$R[Jt]$-module, is either given by ${\mathfrak m}R[Jt]$ or it is empty
$($see \cite[Proposition 2.2]{Vas}$)$. Thus we obtain the following
multiplicity estimate
\[
e_0 = e (\mbox{gr}_I(R)_{+}) \leq e(I/I^2 R[Jt]) + e(S_J(I) \otimes R/I).
\]
Observe that $e(I/I^2 R[Jt])$ is less than or equal to $\lambda(I/I^2)$,
as $I/I^2 R[Jt]$ is the homomorphic image  of the polynomial module
$I/I^2[T_1, \ldots, T_d]$.
On the other hand, $e(S_J(I) \otimes R/I)$ is less than or equal to
$s_0=e(S_J(I))$. Hence $e_0 \leq \lambda(I/I^2) + s_0$. The asserted
inequality now follows from the fact that $s_0 = e_1 - e_0 +
\lambda(R/I)$ $($see \cite[Corollary 3.3]{Vas}$)$.

$(${\it b}$)$
By \cite[Section 22]{Nag}, we only need to show the statement in the
case $d=2$. For $n \gg 0$ the Hilbert--Samuel function of $I$ can be
written as
\begin{equation}\label{In}
\lambda(R/I^n) = e_0 {n+1 \choose 2} - e_1{n \choose 1} + e_2.
\end{equation}
Let $q$ be an integer large enough so that $\lambda(R/I^q)$ is given
by $($\ref{In}$)$ and consider the Hilbert function of $I^q$. For
$n\gg 0 $ one has that
\begin{equation}\label{Inq}
\lambda(R/(I^q)^n) = \widetilde{e}_0 {n+1 \choose 2}-
\widetilde{e}_1{n \choose 1} + \widetilde{e}_2.
\end{equation}
As $\lambda(R/(I^q)^n)=\lambda(R/I^{nq})$, an easy comparison between
$($\ref{In}$)$, with $nq$ in place of $n$, and $($\ref{Inq}$)$ yields
\[
e_0 {nq+1 \choose 2} - e_1{nq \choose 1} + e_2
= \widetilde{e}_0 {n+1 \choose 2}-  \widetilde{e}_1{n \choose 1} +
\widetilde{e}_2
\]
or, equivalently,
\[
\frac{1}{2}q^2e_0n^2 + \left(\frac{1}{2}qe_0-qe_1 \right)n + e_2 =
\frac{1}{2}\widetilde{e}_0n^2 + \left( \frac{1}{2}\widetilde{e}_0 -
\widetilde{e}_1 \right) n +\widetilde{e}_2.
\]
Hence one concludes that
\[
\widetilde{e}_0 = q^2 e_0 \qquad \widetilde{e}_1= qe_1 + \frac{1}{2}
q^2e_0-\frac{1}{2}qe_0 \qquad \widetilde{e}_2 = e_2.
\]
Let $\widetilde{s}_0$ denote the multiplicity of the Hilbert--Samuel
polynomial of the Sally module of $I^q$. By \cite[Corollary 3.3]{Vas}
it satisfies the identity $\widetilde{s}_0 = \widetilde{e}_1 -
\widetilde{e}_0 + \lambda(R/I^q)$. Hence the following calculation goes through
\begin{eqnarray*}
\widetilde{s}_0 & = & \widetilde{e}_1 - \widetilde{e}_0 + \lambda(R/I^q) \\
    & = & \left( qe_1 + \frac{1}{2} q^2e_0-\frac{1}{2}qe_0 \right) -
    q^2e_0 + \left(e_0 {q+1 \choose 2} - e_1{q \choose 1} + e_2
    \right) \\
    & = & e_2.
\end{eqnarray*}
This implies that $e_2=\widetilde{s}_0 \geq 0$, as $\widetilde{s}_0$
is the leading coefficient of a polynomial that measures lengths.
Moreover, $e_2=\widetilde{s}_0=0$ if and only if the Sally module
of $I^q$ is zero, i.e., $I^q$ has reduction number one. \QED

\section{On the depth properties of the associated graded ring
of a class of ${\frak m}$-primary ideals}

In this section we study the depth properties of the associated graded
ring of any ${\frak m}$-primary ideal $I$ in a Cohen--Macaulay local
ring $(R, {\mathfrak m})$ for which there exists a positive integer
$n$ such that $J \cap I^k = JI^{k-1}$ for all $k=1, \ldots, n$ and
$\lambda(I^{n+1}/JI^n) \leq 1$.
In Theorem~\ref{main} we show that the associated graded ring of any
such ideal $I$ has always depth at least $d-1$, where $d$ is the
dimension of the ring $R$, while in Proposition~\ref{cohen-mac} we
single out those ideals whose associated graded ring is
Cohen--Macaulay.

\begin{Proposition}\label{cohen-mac}
Let $(R, \m)$ be a Cohen--Macaulay local ring of dimension $d>0$ and
infinite residue field. Let $I$ be an $\m$-primary ideal of $R$ and
let $J$ be a minimal reduction of $I$ with
$\lambda(I^{n+1}/JI^n)=1$ for some $n \geq 0$. Then the following
conditions are equivalent
\begin{itemize}
\item[$(${\it a}$)$]
$\mbox{\rm gr}_I(R)$ is Cohen--Macaulay;

\item[$(${\it b}$)$]
$J \cap I^k=JI^{k-1}$ for all $k=1, \ldots, n$, $I^{n+1} \not\subset
J$, and $I^{n+2}=JI^{n+1}$.
\end{itemize}
\end{Proposition}
\demo
Suppose that $\mbox{gr}_I(R)$ is Cohen--Macaulay. Then by
\cite[Corollary 2.7]{VV} one has that $J \cap I^k = JI^{k-1}$ for all
$k$. In particular, $J \cap I^{n+1} = JI^n$ and $I^{n+1} \not\subset J$, as
$JI^n\subsetneq I^{n+1}$. Moreover, from the fact that
$\lambda(I^{n+1}/JI^n)=1$ one concludes that $I^{n+2} \subset {\mathfrak m}
I^{n+1} \subset JI^n \subset J$. Hence $I^{n+2}=J \cap
I^{n+2}=JI^{n+1}$.

Conversely, from the short exact sequence
\[
0 \rightarrow J \cap I^{n+1}/ JI^n \longrightarrow I^{n+1}/JI^n
\longrightarrow I^{n+1}/ J \cap I^{n+1} \rightarrow 0
\]
together with the fact that $\lambda(I^{n+1}/JI^n)=1$ and $I^{n+1}/J
\cap I^{n+1} \not= 0$ $($as $I^{n+1} \not\subset J)$ it follows that
$J \cap I^{n+1} = JI^n$. However, $I^{n+2}=JI^{n+1}$ implies that
$J \cap I^k = JI^{k-1}$ for all $k \geq n+2$.
Hence by \cite[Corollary 2.7]{VV} we conclude that the associated
graded ring of $I$ is Cohen--Macaulay. \QED

\medskip

Theorem~\ref{usingtheform} and Theorem~\ref{secondmain}
below describe interesting cases in which condition $(${\it b}$)$ of
Proposition~\ref{cohen-mac} is verified. A lemma is needed first.

\begin{Lemma}\label{makingtheform}
Let $(R, \m)$ be a Cohen--Macaulay local ring of dimension $d>0$ and
infinite residue field. Let $I$ be an $\m$-primary ideal of $R$
and let $J$ be a minimal reduction of $I$. Suppose
$\lambda(I^{n+1}/JI^n)=1$ for some $n\geq 1$ and $I^{n+1} \not\subset
J$. Then
\begin{itemize}
\item[$(${\it a}$)$]
$I^{n+1} \subset (\alpha) +J$ for some $\alpha\in I^{n+1}\setminus
J$;

\item[$(${\it b}$)$]
$V=\widetilde{I}=I+ J\colon I^n/J\colon I^n$ is a finite dimensional
$R/\m$-vector space;

\item[$(${\it c}$)$]
by letting
\[
V^{n+1} \ni (\widetilde{i}_1, \ldots,  \widetilde{i}_{n+1}) \mapsto
f(i_1, \ldots, i_{n+1}) + \m,
\]
where $i_1\cdots i_{n+1} - f(i_1, \ldots, i_{n+1})\alpha \in J$, one
defines a non-degenerate, symmetric, $(n+1)$-linear form $f$ on $V^{n+1}$.
\end{itemize}
If, in addition, the vector space dimension of $V$ is at least $2$
then $\alpha I \subseteq JI^{n+1}$.
\end{Lemma}
\demo
$(${\it a}$)$
By the proof of Proposition~\ref{cohen-mac} we have that $I^{n+1} \cap J =
JI^n$, hence $\lambda(I^{n+1}+J/J)=\lambda(I^{n+1}/JI^n)=1$. Therefore
$I^{n+1}+J/J = (\alpha)+J/J$ for some $\alpha \in I^{n+1}\setminus J$.

$(${\it b}$)$
$V=\widetilde{I}=I+J\colon I^n/J\colon I^n$ is an $R/\m$-vector space since
$I^{n+1}/JI^n$ is.

$(${\it c}$)$
For any given $(\widetilde{i}_1, \ldots,
\widetilde{i}_{n+1})\in V^{n+1}$ it follows that $i_1 \cdots i_{n+1} \in
I^{n+1} \subset (\alpha)+J$. Thus one has
$i_1\cdots i_{n+1}-f(i_1, \ldots, i_{n+1})\alpha\in J$, for some
$f(i_1, \ldots, i_{n+1}) \in R$. By letting
\[
(\widetilde{i}_1, \ldots, \widetilde{i}_{n+1}) \mapsto f(i_1, \ldots,
i_{n+1})+\m
\]
one obtains a well-defined, non-degenerate, symmetric $(n+1)$-linear
form on $V^{n+1}$. We only check that it is well-defined and
non-degenerate, the other properties being trivial.
For the well-definiteness it is enough to show that if for all $t=1,
\ldots, n+1$ one has that $(\widetilde{i}_1, \ldots, \widetilde{i}_t,
\ldots, \widetilde{i}_{n+1})$ and $(\widetilde{i}_1, \ldots,
\widetilde{l}_t, \ldots, \widetilde{i}_{n+1})$ are two representatives
of the same $(n+1)$-tuple of $V^{n+1}$, i.e., $i_t-l_t \in J\colon
I^n$, then $f(i_1, \ldots, i_t, \ldots, i_{n+1}) - f(i_1, \ldots, l_t,
\ldots, i_{n+1}) \in \m$. By assumption we have that
\[
i_1 \ldots i_{t-1}i_{t+1} \ldots i_{n+1}(i_t-l_t) - (f(i_1, \ldots,
i_t, \ldots, i_{n+1}) - f(i_1, \ldots, l_t, \ldots, i_{n+1}))\alpha
\]
is an element in $J$. But $i_1 \ldots i_{t-1}i_{t+1} \ldots
i_{n+1}(i_t-l_t) \in J$, so that
\[
(f(i_1, \ldots, i_t, \ldots, i_{n+1}) - f(i_1,
\ldots, l_t, \ldots, i_{n+1}))\alpha \in J
\]
as well. Therefore, $f(i_1, \ldots, i_t, \ldots,
i_{n+1}) - f(i_1, \ldots, l_t, \ldots, i_{n+1})$ cannot be an
invertible element of $R$, as otherwise this implies $\alpha \in J$.
For the non-degeneracy of the form, suppose that for any $t=1, \ldots,
n+1$ one has $f(i_1, \ldots, i_t, \ldots, i_{n+1})\in \m$
for all $i_j\in I$ with $j\not= t$. By the definition of the form, this
implies that $i_1 \cdots i_t \cdots i_{n+1} \in J$ for all $i_j\in I$ with
$j\not=t$. Hence $i_t \in J \colon I^n$ or, equivalently,
$\widetilde{i}_t=0$.

Finally, if the vector space dimension of $V$ is at least $2$ for any
$c \in I$  one can find an element $\widetilde{d}_2 \in V$ such that
$f(c, d_2,\tratto, \ldots, \tratto) \in \m$.
By the non-degeneracy of the form we can also find $\widetilde{d}_1,
\widetilde{d}_3, \ldots \widetilde{d}_{n+1}$ in $V$ such that $f(d_1,
d_2, d_3, \ldots, d_{n+1})=1$; it follows that $d_1d_2d_3 \cdots
d_{n+1} - \alpha \in I^{n+1} \cap J = JI^n$. Hence $cd_1d_2d_3 \cdots
d_{n+1} - c\alpha \in JI^{n+1}$. On the other hand, $f(c, d_2, d_3,
\ldots, d_{n+1}) \in {\mathfrak m}$ implies $cd_2d_3 \cdots d_{n+1}
\in I^{n+1} \cap J = JI^n$. Hence $d_1(cd_2d_3 \cdots d_{n+1}) \in
JI^{n+1}$, thus yielding $c\alpha \in JI^{n+1}$ as desired. \QED

\medskip

\begin{Theorem}\label{usingtheform}
Let $(R, \m)$ be a Cohen--Macaulay local ring of dimension $d>0$ and
infinite residue field. Let $I$ be an $\m$-primary ideal of $R$
and let $J$ be a minimal reduction of $I$ with
$\lambda(I^{n+1}/JI^n)=1$ for some $n\geq 1$. If the following hold
\begin{itemize}
\item[$(${\it a}$)$]
$J \cap I^k = JI^{k-1}$ for all $k=1, \ldots, n$;

\item[$(${\it b}$)$]
the vector space dimension of $V=I+J\colon I^n/J\colon I^n$ is at least $2$;
\end{itemize}
then the associated graded ring of $I$ is Cohen--Macaulay.
\end{Theorem}
\demo
By Proposition~\ref{cohen-mac} we only need to check that $I^{n+1}
\not\subset J$ and $I^{n+2}=JI^{n+1}$.
If $I^{n+1} \subset J$ then $I \subset J \colon I^n$; this forces the
vector space $V$ to be zero thus contradicting the assumption on its
dimension. Let now $i_1, i_2, \ldots,
i_{n+1}, i_{n+2}$ be $n+2$ arbitrary elements of $I$. If one of them,
say $i_1$, belongs to $J\colon I^n$ then we have that $i_1(i_2 \cdots
i_{n+1}) \in J \cap I^{n+1}=JI^n$ $($by the proof of
Proposition~\ref{cohen-mac}$)$ and $i_1i_2 \cdots i_{n+1}i_{n+2} \in
JI^{n+1}$. Therefore, we may assume that none of the $i_k$'s is in $J
\colon I^n$. The first $n+1$ of them, for example, define a non-zero
element $(\widetilde{i}_1, \ldots, \widetilde{i}_{n+1})$ of
$V^{n+1}$. Making use of the terminology and the results in
Lemma~\ref{makingtheform} we have that $i_1 \cdots i_{n+1} - f(i_1,
\ldots, i_{n+1})\alpha \in J \cap I^{n+1}=JI^n$. Therefore, $i_1
\cdots i_{n+1}i_{n+2} - f(i_1, \ldots, i_{n+1})\alpha i_{n+2} \in
JI^{n+1}$. However, the dimension of $V$ is at least $2$ so that
$\alpha I \subset JI^{n+1}$ by Lemma~\ref{makingtheform}; in particular
$\alpha i_{n+2}\in JI^{n+1}$. Thus $i_1 \cdots i_{n+1}i_{n+2} \in
JI^{n+1}$ as well. \QED

\medskip

\begin{Remark}
{\rm
In the case of Theorem~\ref{usingtheform} with $I=\m$ and $n=1$ $($so
that $e=\mu(\m)-d+2)$, an easy length comparison in the short
exact sequence
\[
0 \rightarrow J\colon \m/J \longrightarrow \m/J \longrightarrow
V = \m/J\colon \m \rightarrow 0
\]
yields the following set of equalities
\begin{eqnarray*}
\mbox{\rm dim}(\m/J\colon \m) & = & \lambda(\m/J\colon \m) =
\lambda(\m/J)-\lambda(J\colon \m/J) \\
& = & e - 1 - \mbox{\rm type}(R) = \mu(\m) - d +1 - \mbox{\rm type}(R).
\end{eqnarray*}
Hence $\mbox{\rm dim}(\m/J\colon \m)\geq 2$ is equivalent to $\mu(\m)-d
> \mbox{\rm type}(R)$. Therefore Theorem~\ref{usingtheform} recovers
the result in \cite{S4}. }
\end{Remark}

\begin{Lemma}\label{utri}
Let $(R, {\mathfrak m})$ be a Gorenstein local ring of dimension $d>0$
and infinite residue field. Let $J$ be a minimal reduction of
${\mathfrak m}$ such that ${\mathfrak m}^{n+1}+J/J$ is the socle of
the Gorenstein ring $R/J$. Then
\[
\lambda({\mathfrak m}/J \colon {\mathfrak m}^n) = \lambda({\mathfrak
m}^n+J/{\mathfrak m}^{n+1}+J).
\]
\end{Lemma}
\demo
As $R$ is Gorenstein one has that $J\colon ({\mathfrak m}^n+J)/J$ is
isomorphic to the canonical module of $R/({\mathfrak m}^n+J)$. Hence
it follows that $\lambda(J \colon {\mathfrak m}^n/J) = \lambda(J\colon
({\mathfrak m}^n+J)/J) = \lambda(R/{\mathfrak m}^n+J) =
\lambda({\mathfrak m}/{\mathfrak m}^n+J)+1$. Thus, an
easy diagram chase together with the fact that ${\mathfrak
m}^{n+1}+J/J$ is the socle of $R/J$ yields
\begin{eqnarray*}
\lambda({\mathfrak m}/J \colon {\mathfrak m}^n) & = &
\lambda({\mathfrak m}^n+J/J) - 1 \\
& = & \lambda({\mathfrak m}^n+J/{\mathfrak m}^{n+1}+J)
+\lambda({\mathfrak m}^{n+1}+J/J) - 1 \\
& = & \lambda({\mathfrak m}^n+J/{\mathfrak m}^{n+1}+J),
\end{eqnarray*}
as claimed. \QED

\begin{Theorem}\label{secondmain}
Let $(R, {\mathfrak m})$ be a Gorenstein local ring of dimension $d>0$
and infinite residue field. Let $J$ be a minimal reduction of
${\mathfrak m}$ with $\lambda({\mathfrak m}^3/J{\mathfrak m}^2)=1$.
Then the associated graded ring of ${\mathfrak m}$ is Cohen--Macaulay.
\end{Theorem}
\demo
We first observe that $J \cap {\mathfrak m}^2 = J{\mathfrak m}$, by the
analytic independence of the generators of $J$, and ${\mathfrak
m}^3 \not\subset J$. Indeed, if ${\mathfrak m}^3 \subset J$ then $R$
has multiplicity $e = n-d+2$, where $n$ denotes the embedding dimension
of $R$. In this case $\lambda({\mathfrak
m}^3/J{\mathfrak m}^2)=0$ by \cite{S3.5}, thus contradicting our
assumption.

By Lemma~\ref{utri}, with $n=2$, we have that $\lambda({\mathfrak m}/ J
\colon {\mathfrak m}^2) = \lambda({\mathfrak m}^2+J/{\mathfrak
m}^3+J)$. By Theorem~\ref{usingtheform} the statement is taken care of if the
previous length is greater than or equal to $2$. Hence we only have to
consider the case in which $\lambda({\mathfrak m}^2+J/{\mathfrak
m}^3+J)=1$. But this condition implies that $R$ has multiplicity
$e=n-d+3$; in this case the Cohen--Macaulayness of ${\rm
gr}_{\mathfrak m}(R)$ follows from \cite[Theorem 1]{S3.75}.
\QED

\bigskip

In Theorem~\ref{secondmain}, the hypothesis of $R$ being Gorenstein
cannot be dropped; moreover, there are examples of Gorenstein rings
with $\mbox{\rm gr}_{\mathfrak m}(R)$ not Gorenstein.

\begin{Example}
{\rm
\begin{itemize}
\item[$(${\it a}$)$]
Let $k$ be a field. The ring $k[\![ t^6, t^7, t^9, t^{17} ]\!]$ is
Cohen--Macaulay, but not Gorenstein, with $\lambda({\mathfrak
m}^3/J{\mathfrak m}^2)=1$, where $J=(t^6)$. In this case, the
associated graded ring of ${\mathfrak m}$ is not Cohen--Macaulay.

\item[$(${\it b}$)$]
Let $k$ be a field. The ring $k[\![ t^5, t^6, t^9 ]\!]$ is Gorenstein
with $\lambda({\mathfrak m}^3/J {\mathfrak m}^2)=1$, where $J = (t^5)$.
By Theorem~\ref{secondmain} the associated graded ring of ${\mathfrak
m}$ is Cohen--Macaulay. However, it is not Gorenstein.
\end{itemize} }
\end{Example}

\medskip

The lemmata below are inspired by and at the same time generalize the
following list of results: \cite[Lemma 2.1 and Lemma 2.2]{Huc1},
\cite[Lemma 1.1, Proposition 1.2, and Corollary 2.3]{RV},
\cite[Lemma 2.1.2]{Vaz}, and \cite[Lemma 2.1, Lemma 2.3, and
Corollary 2.7]{W2}.

\begin{Lemma}\label{lenght}
Let $(R, {\frak m})$ be a Cohen--Macaulay local ring and let $I$ be an
${\frak m}$-primary ideal of $R$. Let $J$ be a minimal reduction of
$I$ such that $\lambda(I^{n+1}/JI^n)=1$. Then either there exists $z
\in I$ such that $I^{t+1} = JI^t + (z^{t+1})$ for all $t \geq n$, or
$I^{n+2}=JI^{n+1}$. In particular, $\lambda(I^{t+1}/JI^t) \leq 1$ for
all $t \geq n  $.
\end{Lemma}
\demo
Let us write $I = (J, z_1, \ldots, z_l)$ for some $l \geq 1$. If there
exists a $k$ such that $z_k^{n+1} \not\in JI^n$ then we may set
$z=z_k$ and we are done. Otherwise, suppose that
$z_i^{n+1} \in JI^n$ for all $i=1, \ldots, l$. Since $JI^n \subsetneq
I^{n+1}$ there exists $f$ such that $I^{n+1} = JI^n + (f)$ of the form
$f = z_1^{p_1}\cdots z_l^{p_l}$ with $p_1 + \cdots + p_l = n+1$.
Choose $i$ such that $p_i>0$ is maximal with respect to the property
that $f \not\in JI^n$. By assumption $p_i < n+1$ so that there exists
$j\not= i$ with $p_j>0$. Note that $(z_1, \ldots, z_l)f \in z_i
I^{n+1} = z_i(JI^n + (f)) \subseteq JI^{n+1} + (z_if) = JI^{n+1} +
(z_j g) =JI^{n+1}$, as $z_if = z_jg$ and $g \in JI^n$ by the choice of $f$.
This implies that
\[
JI^{n+1} \subseteq I^{n+2} = I I^{n+1} = (J, z_1, \ldots, z_l)(JI^n
+(f)) \subseteq JI^{n+1} + (z_1, \ldots, z_l)f = JI^{n+1},
\]
as claimed.

To complete the proof, let us assume that there exists $z \in I$ such
that $I^{n+1} = JI^n + (z^{n+1})$. We will show that $I^{t+1} = JI^t +
(z^{t+1})$ for any $t \geq n$ by inducting on the difference $t-n \geq
0$. If $t-n=0$ there is nothing to prove. Hence by inductive
hypothesis we have
\[
JI^{t+1} + (z^{t+2}) \subseteq I I^{t+1} = I (JI^t +
(z^{t+1})).
\]
Since $I = (J, z_1, \ldots, z_l)$ and  $(z_1, \ldots,
z_l)z^n \subseteq I^{n+1} = JI^n + (z^{n+1})$ one can also write
\begin{eqnarray*}
I(JI^t + (z^{t+1})) & = & J(JI^t + (z^{t+1})) +
(z_1, \ldots, z_l)z^{t+1} \subseteq JI^{t+1} + z^{t+1-n}(z_1, \ldots,
z_l)z^n \\
& \subseteq & JI^{t+1} + z^{t+1-n}(JI^n + (z^{n+1})) = JI^{t+1} +
(z^{t+2}).
\end{eqnarray*}
Hence $I^{(t+1)+1} = JI^{t+1} + (z^{(t+1)+1})$ as requested. The
assertion on  $\lambda(I^{t+1}/JI^t)$ for all $t \geq n$ is now
obvious.
\QED

\bigskip

For the definition and properties of superficial elements$/$sequences
see \cite[Section 22]{Nag}.

\begin{Lemma}\label{lemma-3}
Let $(R, {\frak m})$ be a two dimensional Cohen--Macaulay local ring
and let $I$ be an ${\frak m}$-primary ideal of $R$. Let $J$ be a
minimal reduction of $I$, assume that $J=(x, y)$ where both $x$ and
$y$ are superficial for $I$, and set $r=r_J(I)$ ,
$s=r_{J/(x)}(I/(x))$. If $\lambda(I^{n+1}/JI^n)=1$ and $J \cap I^k =
JI^{k-1}$ for all $k=1, \ldots, n$, then the following statements
hold
\begin{itemize}
\item[$(${\it a}$)$]
$e_1 = \sum\limits_{t=1}^s \lambda(I^t/JI^{t-1})$;

\item[$(${\it b}$)$]
$\mbox{\rm depth}(\mbox{\rm gr}_I(R)) \geq 1$ if and only if $s=r$.
\end{itemize}
\end{Lemma}
\demo
$(${\it a}$)$
As $x$ is a superficial element $e_1=e_1(I)=e_1(I/(x))$. Moreover,
by \cite[Corollary 4.13]{HM} one has that
\begin{eqnarray*}
e_1(I/(x)) \! & = & \!
\sum_{t=1}^s \lambda((I/(x))^t/(J/(x))(I/(x))^{t-1}) = \sum_{t=1}^s
\lambda(I^t+(x)/JI^{t-1}+(x)) \\
\! & = & \! \sum_{t=1}^s \lambda(I^t/JI^{t-1}+((x) \cap I^t)).
\end{eqnarray*}
However, by assumption one has that $(x) \cap I^t \subseteq JI^{t-1}$
for $t=1, \ldots, n$. On the other hand, for $t = n+1, \ldots, s$ it
follows from Lemma~\ref{lenght} that
\[
0 < \lambda((I/(x))^t/(J/(x))(I/(x))^{t-1}) =
\lambda(I^t/JI^{t-1}+ ((x) \cap I^t)) \leq \lambda(I^t/JI^{t-1})
\leq 1,
\]
which implies that $(x) \cap I^t \subseteq JI^{t-1}$ for $t=n+1,
\ldots, s$ as well. This yields the conclusion.

$(${\it b}$)$ The statement follows from \cite[Theorem 3.1]{Huc1}.
\QED

\bigskip

The next theorem contains the third main result of this paper. Its
proof is a simplified version of the one of \cite[Theorem
2.6]{Huc2}, which in turn was inspired by and follows the steps of the
one of \cite[Theorem 2.5]{RV}. The result requires three ingredients:
$(${\it a}$)$ a reduction to the two dimensional case; $(${\it b}$)$
the fact that $e_1$ can be written in two different ways $($one
using the $I$-adic filtration of $I$ and the other using the
filtration given by the Ratliff--Rush closure of the powers of $I$,
see \cite{RR}$)$; and $(${\it c}$)$ a key reduction bound due to
Rossi--Valla.

\begin{Theorem}\label{main}
Let $(R, {\frak m})$ be a Cohen--Macaulay local ring of dimension
$d > 0$ and infinite residue field. Let $I$ be an ${\frak m}$-primary
ideal of $R$ with a minimal reduction $J$ satisfying the following conditions
\begin{itemize}
\item[$(${\it a}$)$]
$J \cap I^k = JI^{k-1}$ for all $k=1, \ldots, n$;

\item[$(${\it b}$)$]
$\lambda(I^{n+1}/JI^n) \leq 1$.
\end{itemize}
Then $\mbox{\rm depth}(\mbox{\rm gr}_I(R)) \geq d-1$.
\end{Theorem}
\demo
By \cite[Lemma 2.2]{HM} the conclusion of the theorem holds in $R$ if
and only if it holds in $R/(x_1, \ldots, x_{d-2})$,
where $x_1, \ldots, x_{d-2} \in J$ is a superficial sequence for $I$:
this is the so called {\it Sally machine}. Moreover, conditions
$(${\it a}$)$ and $(${\it b}$)$ are preserved modulo ${\bf x}=(x_1, \ldots,
x_{d-2})$. Clearly, $\lambda((I/{\bf x})^{n+1}/(J/{\bf x})(I/{\bf
x})^n) \leq 1$ $($see the proof of Lemma~\ref{lemma-3}$)$
so we only need to verify that $(J/{\bf x}) \cap (I/{\bf x})^k =
(J/{\bf x})(I/{\bf x})^{k-1}$ for all $k=1, \ldots, n$. But for that
it will be enough to show that $(J/(x)) \cap (I/(x))^k = JI^{k-1}+(x)/(x)$
holds for any $x \in J$ and for all $k=1, \ldots, n$. Let
$\overline{\imath} = j + ax = i_k + bx$ for some $j \in J$,
$i_k \in I^k$, and $a, b\in R$. Thus
$i_k \in J$ and then, by assumption, $i_k \in J \cap I^k=JI^{k-1}$.
Hence $\overline{\imath} \in JI^{k-1}+(x)/(x)$.

Therefore we may assume $R$ to be two dimensional and $J =
(x, y)$ with $x, y$ superficial elements for $I$. Let
$s=r_{J/(x)}(I/(x))$ and $r=r_J(I)$ as in Lemma~\ref{lemma-3}.
If $r \leq n$ the associated graded ring of $I$ is Cohen--Macaulay
by Valabrega--Valla (see \cite[Corollary 2.7]{VV}). If $r=n+1$ then
the associated graded ring has depth at least $1$ $($or $d-1$
after lifting back$)$ by \cite[Theorem 3.2]{Gu}.
Thus we may assume $r \geq n+2$. The proof will be completed
once we show that $s=r$ $($see Lemma~\ref{lemma-3}$(${\it b}$))$.
By Lemma~\ref{lenght}, there exists $z\in I$ such that $I^{t+1}
= JI^t + (z^{t+1})$ for all $t \geq n$ and $\lambda(I^{t+1}/JI^t) \leq 1$
for all $t \geq n$ $($equality holds if in addition $t<r)$.
The integers
\[
p = \inf \{ k : J\widetilde{I^k} = \widetilde{I^{k+1}} \}
\qquad
q = \inf \{ k : I^{k+1} \subseteq J \widetilde{I^k} \}
\]
satisfy the following inequalities
\[
n \leq q \leq p \leq s.
\]
Indeed, if $q < n$ we have that $I^{q+1} = J\cap I^{q+1}=JI^q$
as $I^{q+1} \subseteq J\widetilde{I^q} \subseteq J$.
But this contradicts the fact that $r \geq n+2$. Hence $n \leq q$.
Since $I^{p+1} \subseteq \widetilde{I^{p+1}} = J\widetilde{I^p}$
it also follows that $q \leq p$. In order to prove the last inequality
notice that $J\widetilde{I^{k-1}} \cap I^k = JI^{k-1}$ for $1 \leq
k \leq n$.
Hence, we obtain the following family of short exact sequences
\[
0 \rightarrow J \widetilde{I^{k-1}}/JI^{k-1}
\stackrel{\varphi_k}{\longrightarrow} \widetilde{I^k}/I^k
\longrightarrow \widetilde{I^k}/J\widetilde{I^{k-1}}+I^k
\rightarrow 0.
\]
Therefore for $k = 2, \ldots, n$ we have that the following expression
\[
\lambda(\widetilde{I^k}/J\widetilde{I^{k-1}}+I^k) =
\lambda(\widetilde{I^k}/I^k) -
\lambda(J\widetilde{I^{k-1}}/JI^{k-1}) =
\lambda(\widetilde{I^k}/J\widetilde{I^{k-1}}) -
\lambda(I^k/JI^{k-1})
\]
is positive. Moreover, $\lambda(\widetilde{I}/I) =
\lambda(\widetilde{I}/J)-\lambda(I/J) \geq 0$, as $I \subseteq
\widetilde{I}$. Consider now the identity
\[
e_1 = \sum_{k=1}^s \lambda(I^k/JI^{k-1}) = \sum_{k \geq 1}
\lambda(\widetilde{I^k}/J\widetilde{I^{k-1}})
\]
that holds by \cite[Corollary 2.10]{Huc1} and
Lemma~\ref{lemma-3}$(${\it a}$)$. We can rewrite the previous formula as
\begin{equation}\label{twicee1}
\sum_{k=1}^n \left(
\lambda(\widetilde{I^k}/J\widetilde{I^{k-1}}) -
\lambda(I^k/JI^{k-1})
\right) = \sum_{k=n+1}^s \lambda(I^k/JI^{k-1}) -
\sum_{k \geq n+1} \lambda(\widetilde{I^k}/J\widetilde{I^{k-1}}).
\end{equation}
>From $($\ref{twicee1}$)$ and Lemma~\ref{lenght} one concludes that
\[
0 \leq s - n - \sum_{k\geq n+1}\lambda(\widetilde{I^k}/J\widetilde{I^{k-1}}).
\]
Hence $s \geq p$. Let $\mu_k$ denote the minimal number of generators
of $\widetilde{I^k}/J\widetilde{I^{k-1}}+I^k$, for each $k \geq 1$.
We have that
\[
\mu_k < \lambda(\widetilde{I^k}/J\widetilde{I^{k-1}}) \qquad \mbox{for
all } k=1, \ldots, q
\]
and also
\[
\mu_k \leq \lambda(\widetilde{I^k}/J\widetilde{I^{k-1}}) -
        \lambda(I^k/JI^{k-1})
\qquad \mbox{for all } k=1, \ldots, n.
\]
By letting $\mu = \mu_1 + \cdots + \mu_q$ it follows from
\cite[Proposition 2.4 and Proof of Theorem 2.5]{RV} $($see also
\cite[Proposition 2.3 and Proof of Theorem 2.6]{Huc2}$)$ that
\[
I^{\mu+q+1} = JI^{\mu+q};
\]
this means that $r \leq \mu + q$. We now show that $\mu+q \leq s$
thus yielding $s=r$. From $($\ref{twicee1}$)$ we have that
\begin{eqnarray*}
\mu & = & (\mu_1 + \cdots + \mu_n) + (\mu_{n+1} + \cdots + \mu_q) \\
    & \leq &
\sum_{k=n+1}^s \lambda(I^k/JI^{k-1}) -
\sum_{k \geq n+1} \lambda(\widetilde{I^k}/J\widetilde{I^{k-1}}) +
(\mu_{n+1} + \cdots + \mu_q) \\
    & = &
s-n + \sum_{k=n+1}^q \left( \mu_k -
    \lambda(\widetilde{I^k}/J\widetilde{I^{k-1}})
\right) - \sum_{k \geq q+1}
\lambda(\widetilde{I^k}/J\widetilde{I^{k-1}}) \\
    & \leq &
s-n - (q-n) - \sum_{k \geq q+1}
\lambda(\widetilde{I^k}/J\widetilde{I^{k-1}})  \\
    & \leq &
s-q
\end{eqnarray*}
or, equivalently, $\mu+q \leq s$.
\QED

\medskip

It is worth pointing out the following consequence of
Theorem~\ref{main} as it describes a situation that
quite frequently occurs in {\it nature}. We note that this result was
previously known only in the case of an ideal $I$ with reduction
number two $($see \cite[Proposition 5.1.4$(${\it a}$)$]{Vas2}$)$.

\begin{Corollary}\label{i2/ji}
Let $(R, {\frak m})$ be a Cohen--Macaulay local ring of dimension
$d > 0$ and infinite residue field. Let $I$ be an ${\frak m}$-primary
ideal of $R$ with a minimal reduction $J$ such that
$\lambda(I^2/JI)=1$. Then the associated graded ring of $I$ has depth
at least $d-1$.
\end{Corollary}

\section{Classes of Examples}

We now describe two situations where the previous results apply.

\subsection{Stretched Cohen--Macaulay rings}

Let $(R, {\frak m})$ be a Cohen--Macaulay local ring
with dimension $d$, infinite residue field,
multiplicity $e$, and embedding dimension $\mu({\frak m})=e
+d-n$ for some $n \geq 1$. The ring $R$ is said to be {\it
stretched} if there exists a minimal reduction $J$ of ${\frak m}$ such
that ${\frak m}^n \not\subset J$.
As $R/J$ is an Artin local ring of length $e$ and embedding
dimension $e-n$ one has that ${\frak m}^{n+1} \subseteq J$.

By combining Theorem~\ref{main} and some results of J. Sally, one
obtains the following

\begin{Proposition}
Let $(R, {\frak m})$ be a Cohen--Macaulay local ring with dimension
$d$, infinite residue field, and embedding dimension $e+d-n$. Let $J$
be a minimal reduction of ${\frak m}$ such that ${\frak m}^n \not\subset J$.
\begin{itemize}
\item[$(${\it a}$)$]
$(${\rm Sally, see \cite[Corollary 2.4]{S3}}$)$
The associated graded ring $\mbox{\rm gr}_{{\frak m}}(R)$ of ${\frak m}$
is Cohen--Macaulay if and only if ${\frak m}^{n+1} = J {\frak m}^n$.

\item[$(${\it b}$)$]
If $\lambda({\frak m}^n / J {\frak m}^{n-1})=1$ then the depth of the
associated graded ring $\mbox{\rm gr}_{{\frak m}}(R)$ of ${\frak m}$
is at least $d-1$.
\end{itemize}
\end{Proposition}
\demo
$(${\it b}$)$
By \cite[Theorem 2.3]{S3} it follows that $J \cap {\mathfrak m}^k = J
{\mathfrak m}^{k-1}$ for all $k =1, \ldots, n$, since ${\mathfrak
m}^{n+1} \subset J {\mathfrak m}^{n-1}$. Hence Theorem~\ref{main}
applies. \QED

\subsection{Ideals arising from graphs}

Let $k$ be a field and let $R$ be the polynomial ring over $k$ in the
$d=2n+1$ variables $x_1, \ldots, x_{n+2}, y_1, \ldots, y_{n-1}$, where
$n\geq 1$. Let
\[
\MM = (x_1, \ldots, x_{n+2}, y_1, \ldots, y_{n-1}) \qquad
\JJ = (x_1^2, \ldots, x_{n+2}^2, y_1^2, \ldots, y_{n-1}^2)
\]
and define a new family of ideals, say $\II_{(3, 2(n-1))}$, as follows
\[
\II_{(3, 2(n-1))} =\left( \JJ , \MM^3, f \right),
\]
where $f$ is the form of degree two given by
\[
f = \sum_{i=1}^{n+1} x_i x_{i+1} + x_{n+2} x_1 + \sum_{i=1}^{n-1}
    x_iy_i.
\]
The ideal $\II_{(3, 2(n-1))}$, or $\II$ for short,
has a simple combinatorial description which arises from the following
graph
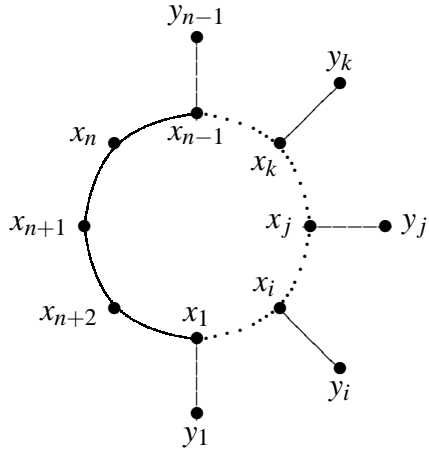
\[
\xy
\POS(-15,0);(0,15)
\curve{(-15,0)&(-14,14)&(0,15)}
\POS(0,15);(15,0)
\curve{~*=<4pt>{.}(0,15)&(14,14)&(15,0)}
\POS(15,0);(0,-15)
\curve{~*=<4pt>{.}(15,0)&(14,-14)&(0,-15)}
\POS(-15,0);(0,-15)
\curve{(-15,0)&(-14,-14)&(0,-15)}
\POS(-15,0)*+{\bullet}
\POS(0,15)*+{\bullet}
\POS(0,-15)*+{\bullet}
\POS(-11,11)*+{\bullet}
\POS(-11,-11)*+{\bullet}
\POS(-15,12)*+{x_n}
\POS(-21,0)*+{x_{n+1}}
\POS(-17,-12)*+{x_{n+2}}
\POS(0,-12)*+{x_1}
\POS(0,12)*+{x_{n-1}}
\POS(0,15);(0,25)
\curve{}
\POS(0,-15);(0,-25)
\curve{}
\POS(0,25)*+{\bullet}
\POS(0,-25)*+{\bullet}
\POS(0,-28)*+{y_1}
\POS(0,28)*+{y_{n-1}}
\POS(11,11)*+{\bullet}
\POS(15,0)*+{\bullet}
\POS(11,-11)*+{\bullet}
\POS(9,8)*+{x_k}
\POS(11,0)*+{x_j}
\POS(9,-8)*+{x_i}
\POS(19,22)*+{y_k}
\POS(29,0)*+{y_j}
\POS(19,-22)*+{y_i}
\POS(19,19)*+{\bullet}
\POS(25,0)*+{\bullet}
\POS(19,-19)*+{\bullet}
\POS(11,11);(19,19)
\curve{}
\POS(15,0);(25,0)
\curve{}
\POS(11,-11);(19,-19)
\curve{}
\endxy
\]
More precisely, the element $f$ is obtained by adding together all the
products $x_ix_{i+1}$, where the pair $(x_i, x_{i+1})$ consists of
adjacent vertices on the cycle, and all the products $x_iy_i$, where
the pair $(x_i, y_i)$ is a {\it whisker}.

\begin{Example}{\rm
If $n=1$ then $2(n-1)=0$ and the ideal $\II=\II_{(3,0)}$ corresponds
to the case of a cycle with three vertices and no whisker.
More precisely
\[
\II = (x_1^2, x_2^2, x_3^2, f),
\]
where $f= x_1x_2+x_2x_3+x_3x_1$.
It turns out that $\lambda(R/\II)=6$, $\lambda(R/\JJ)=8$,
$\lambda(\II/\JJ)=2$, $\lambda(\II^2/\JJ\II)=1$, $\II^2 \subset
\JJ$ so that $\JJ \cap \II^2 = \II^2 \not= \JJ \II$. Hence the depth
of the associated graded ring of $\II$ is exactly $d-1=3-1=2$. }
\end{Example}

\begin{Example}{\rm
If $n=2$ then $2(n-1)=2$ and the ideal $\II=\II_{(3,2)}$ corresponds
to the case of a cycle with four vertices and one whisker.
More precisely
\[
\II = (x_1^2, x_2^2, x_3^2, x_4^2, y_1^2,
x_1x_3x_4, x_2x_3x_4, x_1x_4y_1, x_3x_4y_1, x_2x_4y_1, f),
\]
where $f = x_1x_2+x_2x_3+x_3x_4+x_4x_1+x_1y_1$. This example appears
in \cite[Example 2.13]{Huc2}.
It turns out that $\lambda(R/\II)=15$, $\lambda(R/\JJ)=32$,
$\lambda(\II/\JJ)=17$, $\lambda(\II^2/\JJ\II)=2$, $\JJ \cap \II^2
= \JJ \II$, $\lambda(\II^3/\JJ \II^2)=1$, $\II^3 \subset
\JJ$ so that $\JJ \cap \II^3 = \II^3 \not= \JJ \II^2$. Hence the depth
of the associated graded ring of $\II$ is exactly $d-1=5-1=4$. }
\end{Example}

\begin{Example}{\rm
If $n=3$ then $2(n-1)=4$ and the ideal $\II=\II_{(3,4)}$ corresponds
to the case of a cycle with five vertices and two whiskers.
More precisely
\[
\II = (x_1^2, x_2^2, x_3^2, x_4^2, x_5^2, y_1^2, y_2^2, \MM^3,
       f),
\]
where $f = x_1x_2+x_2x_3+x_3x_4+x_4x_5+x_5x_1+x_1y_1+x_2y_2$.
It turns out that $\lambda(R/\II)=28$,
$\lambda(R/\JJ)=128$, $\lambda(\II/\JJ)=100$,
$\lambda(\II^2/\JJ\II)=30$, $\JJ \cap \II^2 = \JJ \II$,
$\lambda(\II^3/\JJ \II^2)=2$, $\JJ \cap \II^3 =
\JJ \II^2$, and $\lambda(\II^4/\JJ\II^3)=1$, $\II^4 \subset \JJ$
so that $\JJ \cap \II^4 = \II^4 \not= \JJ\II^3$. Hence the depth
of the associated graded ring of $\II$ is exactly $d-1=7-1=6$. }
\end{Example}

\begin{Question}{\rm
It is natural then to ask whether or not in general the ideals $\II =
\II_{(3, 2(n-1))}$
and $\JJ$ satisfy the following properties: $(${\it a}$)$
$\JJ \cap \II^k = \JJ \II^{k-1}$ for all $k=1 ,\ldots, n$;
$(${\it b}$)$ $\lambda(\II^{n+1}/\JJ \II^n)=1$;
$(${\it c}$)$ $\II^{n+1} \subset \JJ$. }
\end{Question}

\subsection*{Note added in proof}

After this paper was completed, the authors learned that both
M.E. Rossi and J. Elias independently wrote articles that partially
overlap with ours.

\end{document}